\documentclass[12pt]{article}
\usepackage{amsmath,dsfont}
\allowdisplaybreaks[4]

\usepackage[all]{xy}
\usepackage{amssymb}
\usepackage{amsthm}
\usepackage{hyperref}
\hypersetup{colorlinks=true,linkcolor=blue,citecolor=red}
\usepackage{amsmath}
\usepackage{amscd}%,enumitem
\usepackage{verbatim}
\usepackage{eurosym}
\usepackage{float}
\usepackage{color}
\usepackage{dcolumn}
\usepackage[mathscr]{eucal}
\usepackage[all]{xy}
\usepackage{hyperref}
\usepackage{mathrsfs}
\usepackage{amsmath}
\usepackage{amssymb}
\usepackage{amsfonts,ifpdf}
\usepackage{graphicx}
\usepackage{times}
\usepackage{float}
\usepackage{epstopdf}
\usepackage{cite}
\usepackage{youngtab}
\usepackage{ytableau}
\ytableausetup{mathmode, boxsize=0.9em}

\newtheorem{theo}{Theorem}[section]

\newtheorem{lem}[theo]{Lemma}

\newtheorem{coro}[theo]{Corollary}
\newtheorem{con}[theo]{Conjecture}

\newtheorem{pf}{proof}[section]
\theoremstyle{remark}

\makeatletter \@addtoreset{equation}{section} \makeatother
\makeindex 
\def\qed{\hfill \rule{4pt}{7pt}}
\def\pf{\vskip 0.2cm {\noindent \bf Proof.}\quad}
\begin{document}

\begin{center}

 {\Large \bf The $q$-Log-Concavity and Unimodality of\\[8pt]   $q$-Kaplansky Numbers}

\end{center}

\begin{center}
 {Kathy Q. Ji}   \vskip 2mm

   Center for Applied Mathematics,\\[2pt]
Tianjin University\\[2pt]
 Tianjin 300072, P.R. China\\[3pt]

   \vskip 2mm

 kathyji@tju.edu.cn
\end{center}

\vskip 6mm \noindent {\bf Abstract.} $q$-Kaplansky  numbers were considered by Chen and Rota. We find that    $q$-Kaplansky  numbers  are connected to the symmetric differences of   Gaussian polynomials introduced by  Reiner and Stanton. Based on the work of Reiner and Stanton, we establish the  unimodality of   $q$-Kaplansky numbers.   We also show that   $q$-Kaplansky  numbers are the generating functions for  the inversion number and the major index of two special kinds of $(0,1)$-sequences. Furthermore, we  show that   $q$-Kaplansky numbers are  strongly $q$-log-concave.

\noindent
{\bf Keywords:} Inversion number, major index, $q$-log-concavity, unimodality,  $q$-Catalan numbers, Foata's fundamental bijection, integer partitions

\noindent
{\bf AMS Classification:} 05A17, 05A19, 05A20, 05A30

\section{Introduction}

The main objective of this paper is to give two combinatorial interpretations of    $q$-Kaplansky numbers introduced by Chen and Rota \cite{Chen-Rota-1992} and to  establish some properties of  $q$-Kaplansky numbers. Recall that the Kaplansky  number $K(n,m)$ is defined by
\[ K(n,m)={n \over n-m} {n-m \choose m},\]
for $n\geq 2m\geq 0$.  The combinatorial interpretation of $ K(n,m)$ was first given by Kaplansky \cite{Kaplansky-1944},  so we call $ K(n,m)$ the Kaplansky number. Kaplansky found that $ K(n,m)$  counts the number of ways of choosing $m$ nonadjacent elements arranged on a cycle, which can also be interpreted as the number of dissections of type $1^{n-2k}2^k$ of an $n$-cycle given by Chen, Lih and Yeh  \cite{Chen-Lih-Yeh-1995}.    Kaplansky numbers appear in many
 classical polynomials, such as   Chebyshev polynomials of the first kind \cite{Lidl-Mullen-Turnwald-1993, Mason-Handscomb-2003} and  Lucas polynomials \cite{Koshy-2019}.

  $q$-Kaplansky  numbers were introduced by Chen and Rota \cite{Chen-Rota-1992}. For convenience, we adopt the following definition: For $n\geq 1$ and $0\leq m\leq n$,
\begin{equation}\label{defi-Kap}
K_q(n,m)=\frac{1-q^{n+m}}{1-q^{n}}{n \brack m},
\end{equation}
where ${n \brack m} $ is the Gaussian polynomial, also  called the $q$-binomial coefficient, as given by
\[{n \brack m}=\frac{(1-q^n)(1-q^{n-1})\cdots (1-q^{n-m+1})}{(1-q^m)(1-q^{m-1})\cdots (1-q)}.\]
 By the symmetric property of the Gaussian polynomial, it is not hard to show that $K_q({n,m})$ is a symmetric polynomial of degree $m(n-m)+m$ with nonnegative coefficients.

The first result of this paper is to give two combinatorial interpretations of   $q$-Kaplansky numbers.
Let  $w=w_1w_2\cdots w_n$ be a $(0,1)$-sequence of length $n$, the number of inversions of $w$, denoted ${\rm{inv}}(w)$, is the number of pairs $(w_i,w_j)$ such that $i<j$ and $w_i>w_j$, and the major index of $w$, denoted  ${\rm{maj}}(w)$, is  the sum of indices $i<n$ such that $w_i>w_{i+1}$. For example, for $w=10010110$, we have ${\rm{inv}}(w)=8$ and ${\rm{maj}}(w)=1+4+7=12$.

It can be shown that   $q$-Kaplansky numbers are related to two sets  $\mathcal{K}(m,n-m+1)$  and  $ \mathcal{\overline{K}}({m,n-m+1})$  of $(0,1)$-sequences. More precisely, for $n\geq m\geq 0$, let $\mathcal{K}(m,n-m+1)$ denote the set of $(0,1)$-sequences $w=w_1w_2\cdots w_{n+1}$ of length $n+1$ consisting of $m$ copies of 1's and $n-m+1$ copies of 0's such that if $w_{n+1}=1$, then $w_{1}=0$. For $n\geq m\geq 0$, let $ \mathcal{\overline{K}}({m,n-m+1})$ denote  the set of $(0,1)$-sequences $w=w_1w_2\cdots w_{n+1}$ of length $n+1$  consisting of $m$ copies of 1's and $n-m+1$ copies of 0's such that if  $w_{n+1}=1$ and $t:=\max\{i:w_i=0\}$,  then $t=1$ or  $w_{t-1}=0$ when $t\geq 2$.  We have the following combinatorial interpretations.

\begin{theo}\label{main-CombInt} For $n\geq m\geq 0$,
\begin{eqnarray}
K_q({n,m}) &=&\sum_{w \in \mathcal{K}(m,n-m+1)}q^{{{\rm{inv}}}(w)}\label{main-CombInt-e1}\\[5pt]
 &=& \sum_{w \in \mathcal{\overline{K}}({m,n-m+1})} q^{{\rm{maj}}(w)}. \label{main-CombInt-e2}
\end{eqnarray}
\end{theo}

The second result of this paper is to establish  the strong $q$-log-concavity of $K_q(n,m)$. Recall that a sequence
of polynomials $\left(f_n(q)\right)_{n\geq0}$ over the field of real numbers is called $q$-log-concave if the
difference
\[f_m(q)^2-f_{m+1}(q)f_{m-1}(q)\]
has nonnegative coefficients as a polynomial in $q$ for all $m\geq 1$. Sagan \cite{Sagan-1992a} also introduced the notion of the strong $q$-log-concavity. We say that a sequence of polynomials
$\left(f_n(q)\right)_{n\geq0}$  is strongly $q$-log-concave if
\[f_n(q)f_m(q)-f_{n-1}(q)f_{m+1}(q)\]
has nonnegative coefficients  as a polynomial in $q$ for any $m\geq n\geq 1$.

It is known that  $q$-analogues of many well-known combinatorial numbers are
strongly $q$-log-concave. Butler \cite{Butler-1990} and Krattenthaler \cite{Krattenthaler-1989} proved the strong $q$-log-concavity of
  $q$-binomial coefficients, respectively. Leroux \cite{Leroux-1990} and Sagan \cite{Sagan-1992a} studied the strong $q$-log-concavity of $q$-Stirling numbers of the first kind and the second kind. Chen, Wang and Yang \cite{Chen-Wang-Yang-2010} have shown that  $q$-Narayana numbers are strongly $q$-log-concave.

We obtain the following   result which  implies that  $q$-Kaplansky numbers are   strongly $q$-log-concave.

\begin{theo}\label{main-q-log-cancave} For $1\leq m\leq l<n$ and $0\leq r\leq 2l-2m+2$,
\begin{equation}\label{kqtt}
K_q(n,m)K_q(n,l)
-q^{r}K_q(n,m-1)K_q(n,l+1)
\end{equation}
has nonnegative coefficients as a polynomial in $q$.
\end{theo}

\begin{coro}\label{main-q-log-cancave-c} Given a positive integer $n$, the sequence  $\left(K_q(n,m)\right)_{0\leq m\leq n}$ is strongly $q$-log-concave.
\end{coro}

It is easy to check that the degree of $K_q(n,m)K_q(n,l)$ exceeds the degree of $K_q(n,m-1)K_q(n,l+1)$ by $2l-2m+2$, so if the difference \eqref{kqtt} of these two polynomials has nonnegative coefficients, then $r\leq 2l-2m+2$.

To conclude the introduction, let us say a few words about the unimodiality of $q$-Kaplansky numbers. We find that   $q$-Kaplansky numbers are connected to the following symmetric differences  of  Gaussian polynomials introduced by  Reiner and Stanton \cite{Reiner-Stanton-1998}.
\begin{equation}
 F_{n,m}(q)={n+m\brack m}-q^n {n+m-2\brack m-2}.
\end{equation}

The following theorem is due to  Reiner and Stanton \cite{Reiner-Stanton-1998}.

\begin{theo}[Reiner-Stanton] \label{Reiner-Stanton} When $m\geq 2$ and $n$ is even, the polynomial $F_{n,m}(q)$ is symmetric and unimodal.
\end{theo}
Recently, Chen and Jia \cite{Chen-Jia} provided a simple proof of the unimodality of $F_{n,m}(q)$ by using semi-invariants. According to the following recursions of Gaussian polynomials \cite[p.35,Theorem 3.2 (3.3)]{Andrews-1976},
\begin{eqnarray}\label{rec-a}
{n \brack m}&=&{n-1 \brack m-1}+q^m {n-1 \brack m},  \\[8pt]
{n-1 \brack m}&=&{n \brack m}-q^{n-m}{n-1 \brack m-1},   \label{rec-b}
\end{eqnarray}
we  find that
\begin{eqnarray}\label{qkap-ReiSta}
F_{n,m}(q)&=&{n+m\brack m}-q^n {n+m-2\brack m-2} \nonumber\\[5pt]
&\overset{\eqref{rec-a}}{=}& {n+m-1\brack m-1}-q^n {n+m-2\brack m-2}+q^m{n+m-1\brack m} \nonumber \\[5pt]
&\overset{\eqref{rec-b}}{=}& {n+m-2\brack m-1}+q^m{n+m-1\brack m} \nonumber \\[8pt]
&=& \frac{1-q^{n+2m-1}}{1-q^{n+m-1}}{n+m-1\brack m} \nonumber \\[8pt]
&=& K_q({n+m-1,m}).
\end{eqnarray}

 Combining Theorem \ref{Reiner-Stanton} and \eqref{qkap-ReiSta}, we have the following result.

 \begin{theo}\label{Unimodal} When $n\geq m\geq 2$ and $n-m$ is odd,  the $q$-Kaplansky number $K_q({n,m})$ is symmetric and unimodal.
\end{theo}

It should be noted that $K_q({n,m})$ is not always unimodal for any $n\geq m\geq 2$. For example, \[K_q(6,2)=1+q+2q^2+2q^3+3q^4+2q^5+3q^6+2q^7+2q^8+q^9+q^{10}\]
 is not unimodal.

 $q$-Kaplansky numbers are also related to  $q$-Catalan polynomials $C_n(q)$, defined by
 \begin{equation}\label{q-catalan}
C_n(q)=\frac{1-q}{1-q^{n+1}}{2n \brack n}=\frac{1-q}{1-q^{2n+1}}{2n+1 \brack n}.
\end{equation}
  It is well-known that $C_n(q)$ is a polynomial in $q$ with non-negative coefficients \cite{Furlinger-Hofbauer-1968}.
 Combining \eqref{defi-Kap} and \eqref{q-catalan}, it is readily seen that
  \[(1-q) {K}_q({{2n+1}, n})=(1-q^{3n+1})C_n(q).\]

  Hence, by Theorem \ref{Unimodal}, we obtain the following result.

  \begin{theo} When $n$ is even, the polynomial $\frac{1-q^{3n+1}}{1-q}C_n(q)$ is symmetric and unimodal.
\end{theo}

Finally, we would like to state a result of Stanley \cite[p.523]{Stanley-1989}  about the  unimodality of the $q$-Catalan polynomials   and two conjectures on the  unimodality of the $q$-Catalan polynomials due to   Chen, Wang and Wang \cite{Chen-Wang-Wang-2008} and Xin and Zhong \cite[Conjecture 1.2]{Xin-Zhong-2020}, respectively.  Apparently,  Conjecture \ref{conj-chen-wang-wang} implies Conjecture \ref{conj-xin-zhong}  when $n\geq 16.$

 \begin{theo}[Stanley]For $n\geq 1$, the polynomial $\frac{1+q}{1+q^n}C_n(q)$ is symmetric and unimodal.
 \end{theo}

  \begin{con}[Chen, Wang and Wang]\label{conj-chen-wang-wang}
For $n\ge 16$, the $q$-Catalan polynomial $C_n(q)$ is  unimodal.
\end{con}

\begin{con}[Xin and Zhong]\label{conj-xin-zhong}
For $n\ge 1$, the  polynomial $(1+q)C_n(q)$ is unimodal.
\end{con}

\section{Proof of Theorem \ref{main-CombInt}}
To prove Theorem \ref{main-CombInt}, we first recall a result due to MacMahon \cite{MacMahon-1960}. For $n\geq m\geq 0$, let $\mathcal{M}(m,n-m)$ be the set of $(0,1)$-sequences of length $n$ consisting of $m$ copies of 1's and $n-m$ copies of 0's.  The following well-known result is due to MacMahon (see \cite[Chapter 3.4]{Andrews-1976}).

\begin{theo}[MacMahon]\label{Mac} For $n\geq m\geq 0,$
\begin{eqnarray}\label{Mac-e-1}
{n \brack m}&=&\sum_{w\in \mathcal{M}(m,n-m)} q^{{{\rm{ inv}}}(w)}\\[5pt]
&=&\sum_{w \in \mathcal{M}(m,n-m)} q^{{\rm {maj}}(w)}. \label{Mac-e-2}
\end{eqnarray}
\end{theo}
Foata's fundamental bijection  \cite{Foata-1968} can be used to establish the equivalence of \eqref{Mac-e-1} and \eqref{Mac-e-2}. There are several ways to describe Foata's fundamental bijection, see, for example, Foata \cite{Foata-1968}, Haglund \cite[p.2]{Haglund-2008} and Sagan and Savage \cite{Sagan-Savage-2012}. Here we give a description due to Sagan and Savage \cite{Sagan-Savage-2012}.

\noindent{\it Proof of the equivalence between \eqref{Mac-e-1} and \eqref{Mac-e-2}}: Let  $w=w_1w_2\cdots w_{n}\in \mathcal{M}(m,n-m)$. We aim to construct a $(0,1)$-sequence $\widetilde{w}=\phi(w)=\widetilde{w}_1\widetilde{w}_2\cdots \widetilde{w}_{n}$ in  $\mathcal{M}(m,n-m)$  such that ${\rm{inv}}(\widetilde{w})={\rm{maj}}({w}).$

Let $w$ be  a $(0,1)$-sequence with $d$ descents, so that we can write
\begin{equation}
w=0^{m_0}1^{n_0}0^{m_1}1^{n_1}0^{m_2}\cdots 1^{n_{d-1}}0^{m_d} 1^{n_d},
\end{equation}
where $m_0\geq 0$ and $m_i\geq 1$ for $1\leq i\leq d$,  $n_i\geq 1$ for $0\leq i\leq d-1$ and $n_d
\geq 0$.

Define
\begin{equation}
\widetilde{w}=\phi(w)=0^{m_d-1}10^{m_{d-1}-1}1\cdots 0^{m_1-1}10^{m_0}1^{n_0-1}01^{n_1-1}\cdots 01^{n_{d-1}-1}01^{n_d}.
\end{equation}
 It has been shown in  \cite{Sagan-Savage-2012}  that ${\rm{inv}}(\widetilde{w})={\rm{maj}}({w})$.

 The inverse map $\phi^{-1}$ of $\phi$ can be described  recursively. Let $\widetilde{w}\in \mathcal{M}(m,n-m)$, we may write $\widetilde{w}=0^{a}1u01^b$ for $a,b\geq 0$, define
\begin{equation}\label{Foata-r-e}
{w}=\phi^{-1}(\widetilde{w})=\phi^{-1}(u)10^{a+1}1^b.
\end{equation}
 It has been proved in  \cite{Sagan-Savage-2012}  that $\phi^{-1}(\phi({w}))={w}$ and $\phi(\phi^{-1}(\widetilde{w}))=\widetilde{w}$. Furthermore,  ${\rm{inv}}(\widetilde{w})={\rm{maj}}({w}).$   Hence the map $\phi$ is a bijection.  This completes the proof of  the equivalence of \eqref{Mac-e-1} and \eqref{Mac-e-2}. \qed

For $n\geq m\geq 0$, let ${\mathcal{M}}_0({m,n-m+1})$ be the set of $(0,1)$-sequences $w=w_1w_2\cdots w_{n+1}$ of length $n+1$ consisting of $m$ copies of 1's and $n-m+1$ copies of 0's such that $w_{n+1}=0$. We have the following result.

\begin{lem}\label{lemma} For $n\geq m\geq 0,$
\begin{eqnarray}\label{lemma-e-1}
q^m{n \brack m}&=&\sum_{w \in {\mathcal{M}}_0({m,n-m+1})} q^{{\rm{inv}}(w)}\\[5pt]
&=&\sum_{w \in {\mathcal{M}}_0({m,n-m+1})} q^{{\rm{maj}}(w)}. \label{lemma-e-2}
\end{eqnarray}
\end{lem}

\pf By Theorem \ref{Mac}, we see that
\[{n \brack m}=\sum_{w \in \mathcal{M}(m,n-m)} q^{{\rm{inv}}(w)}.\]
To prove \eqref{lemma-e-1}, it suffices to show that
\begin{equation}\label{lemma-e-1-t}
\sum_{w \in \mathcal{M}(m,n-m)} q^{{{\rm{inv}}}(w)+m}=\sum_{w \in {\mathcal{M}}_0({m,n-m+1})} q^{{\rm{inv}}(w)}.
\end{equation}
We construct a bijection $\psi$ between the set $\mathcal{M}(m,n-m)$ and the set ${\mathcal{M}}_0({m,n-m+1})$ such that for $w\in \mathcal{M}(m,n-m)$ and $\psi(w) \in {\mathcal{M}}_0({m,n-m+1})$, we have
\begin{equation}\label{invaa}
{{\rm{inv}}}(w)+m={\rm{inv}}(\psi(w)).
\end{equation}
 Let $w={w}_1{w}_2\cdots {w}_{n}$. Define
 \[\psi(w)={w}_1{w}_2\cdots {w}_{n}0.\]
 It is clear that $\psi(w) \in {\mathcal{M}}_0({m,n-m+1})$ and \eqref{invaa} holds. Furthermore, it is easy to see that $\psi$ is reversible. Hence   we have\eqref{lemma-e-1-t}.

We proceed to show that  \eqref{lemma-e-1} and \eqref{lemma-e-2} are equivalent by using Foata's fundamental bijection $\phi$.
Let $w={w}_1{w}_2\cdots {w}_{n+1}$ be in ${\mathcal{M}}_0({m,n-m+1})$, by definition, we see that ${w}_{n+1}=0$. Define
\[\widetilde{w}=\phi^{-1}({w})=\widetilde{w}_1\widetilde{w}_2\cdots \widetilde{w}_{n+1},\]
 where $\phi^{-1}$ is defined in \eqref{Foata-r-e}. By \eqref{Foata-r-e}, we see that $\widetilde{w}_{n+1}=0$ since ${w}_{n+1}=0$. Hence $\widetilde{w} \in {\mathcal{M}}_0({m,n-m+1})$. Furthermore $\phi^{-1}$ is reversible and ${\rm{inv}}({w})={\rm{maj}}(\widetilde{w}).$ It follows  \eqref{lemma-e-1} and \eqref{lemma-e-2} are equivalent, and so \eqref{lemma-e-2} is valid. \qed

 For $n\geq m\geq 1$, let ${\mathcal{M}}_1({m,n-m+1})$  be the set of $(0,1)$-sequences $w=w_1w_2\cdots w_{n+1}$ of length $n+1$ consisting of $m$ copies of 1's and $n-m+1$ copies of 0's such that $w_1=0$ and $w_{n+1}=1$. For $n\geq m\geq 1$, let $\overline{\mathcal{M}}_1({m,n-m+1})$  be the set of $(0,1)$-sequences $w=w_1w_2\cdots w_{n+1}$ of length $n+1$ consisting of $m$ copies of 1's and $n-m+1$ copies of 0's such that  $w_{n+1}=1$ , and if $t:=\max\{i:w_i=0\}$,  then $t=1$ or  $w_{t-1}=0$ when $t\geq 2$.
 To wit, for  $w \in \overline{\mathcal{M}}_1({m,n-m+1})$, if $m\geq 1$ and $n>m$, then  $w$ can be written as $u001^{n+1-t}$, where $2\leq t\leq n$ and $u \in \mathcal{M}({m+t-n-1,n-m-1})$, and if $m\geq 1$ and $n=m$, then $w$ can be written as $01^{m}$.

 \begin{lem}\label{lemma2} For $n\geq m\geq 1$,
\begin{eqnarray}\label{lemma2-e-1}
 {n-1 \brack m-1}&=& \sum_{w \in \mathcal{\mathcal{M}}_1({m,n-m+1})} q^{{\rm{inv}}(w)}\\[5pt]
 &=& \sum_{w \in \overline{\mathcal{M}}_1({m,n-m+1})} q^{{\rm{maj}}(w)}. \label{lemma2-e-2}
\end{eqnarray}
 \end{lem}
 \pf By Theorem \ref{Mac}, we see that
\[{n-1 \brack m-1}=\sum_{w \in \mathcal{M}({m-1,n-m})} q^{{\rm{inv}}(w)}.\]
To prove \eqref{lemma2-e-1},  it suffices to show that
\begin{equation}\label{lemma2-e-1-ta}
\sum_{w \in \mathcal{M}(m-1,n-m)} q^{{\rm{inv}}(w)}=\sum_{w \in \mathcal{\mathcal{M}}_1({m,n-m+1})} q^{{\rm{inv}}(w)}.
\end{equation}

We now construct a bijection $\varphi$ between the set $\mathcal{M}({m-1,n-m})$ and the set \break $\mathcal{\mathcal{M}}_1({m,n-m+1})$ such that for $w \in \mathcal{M}({m-1,n-m})$ and $\varphi(w) \in \mathcal{\mathcal{M}}_1({m,n-m+1})$, we have
\begin{equation}\label{ttaa}
{\rm{inv}}(w)={\rm{inv}}(\varphi(w)).
\end{equation}
 Let $w={w}_1{w}_2\cdots {w}_{n-1}$. Define
 \[\varphi(w)=0{w}_1{w}_2\cdots {w}_{n-1}1.\]
  It is clear that $\varphi(w) \in \mathcal{\mathcal{M}}_1({m,n-m+1})$ and \eqref{ttaa} holds. Furthermore,  $\psi$ is reversible. Hence   we have \eqref{lemma2-e-1-ta}.

We proceed to show that \eqref{lemma2-e-2} holds. By \eqref{Mac-e-2}, it suffices to show that
\begin{equation}\label{lemma2-e-1-t}
\sum_{w \in \mathcal{M}(m-1,n-m)} q^{{\rm{maj}}(w)}=\sum_{w \in \overline{\mathcal{M}}_1({m,n-m+1})} q^{{\rm{maj}}(w)}.
\end{equation}
We now construct a bijection $\tau$ between the set  $\mathcal{M}({m-1,n-m})$ and the set $\overline{\mathcal{M}}_1(m,$ $n-m+1)$  such that for $w \in \mathcal{M}({m-1,n-m})$ and $\tau(w) \in \overline{\mathcal{M}}_1({m,n-m+1})$, we have
\begin{equation}\label{ttcc}
{\rm{maj}}(w)={\rm{maj}}(\tau(w)).
\end{equation}
 Let $w={w}_1{w}_2\cdots {w}_{n-1} \in \mathcal{M}({m-1,n-m})$. If $n=m$, then $w=1^{m-1}$, and so define $\tau(w)=01^m$.
 If $n>m$, then let $t=\max\{i:w_i=0\}$, obviously, $t\geq 1$.  In this case, we may write $w={w}_1{w}_2\cdots {w}_{t-1}01^{n-t-1} $. Define \[\widetilde{w}=\tau(w)=\widetilde{w}_1\widetilde{w}_2\cdots \widetilde{w}_{n+1}\]
  as follows: set $\widetilde{w}_{n+1}=1$, and set $\widetilde{w}_j={w}_j$ for $1\leq j\leq t$, $\widetilde{w}_{t+1}=0$, and set $\widetilde{w}_{j+1}=w_j=1$ for $t+1\leq j\leq n-1$.

  From the above construction, it is easy to see that $\widetilde{w} \in \overline{\mathcal{M}}_1({m,n-m+1})$ and \eqref{ttcc} holds. Furthermore, it can be checked that this construction is reversible, so \eqref{lemma2-e-1-t} is valid. \qed

We are now in a position to give a proof of Theorem \ref{main-CombInt} based on Lemma \ref{lemma} and Lemma \ref{lemma2}.

\noindent{\it Proof of Theorem \ref{main-CombInt}:} By the definition of $\mathcal{K}(m,n-m+1)$, we see that
\[\mathcal{K}(m,n-m+1)=\mathcal{M}_0(m,n-m+1)\cup \mathcal{\mathcal{M}}_1(m,n-m+1).\]
Combining   \eqref{lemma-e-1} and \eqref{lemma2-e-1}, we derive that for $n\geq m\geq 1$,
\begin{eqnarray*}
\sum_{w \in \mathcal{K}({m,n-m+1})}q^{{{\rm{inv}}}(w)}&=&\sum_{w \in \mathcal{M}_0(m,n-m+1)}q^{{{\rm{inv}}}(w)}+\sum_{w \in \mathcal{M}_1(m,n-m+1)}q^{{{\rm{inv}}}(w)}\\[5pt]
&=&q^m{n\brack m}+{n-1\brack m-1}\\[5pt]
&=&\frac{1-q^{n+m}}{1-q^n}{n\brack m}\\[5pt]
&=&{K}_q({n,m}).
\end{eqnarray*}

Similarly, by definition,  we see that
\[\overline{\mathcal{K}}(m,n-m+1)=\mathcal{M}_0(m,n-m+1)\cup \overline{\mathcal{M}}_1(m,n-m+1).\]
By   \eqref{lemma-e-2} and \eqref{lemma2-e-2}, we find that $n\geq m\geq 1$,
\begin{eqnarray*}
\sum_{w \in \overline{\mathcal{K}}({m,n-m+1})}q^{{{\rm{maj}}}(w)}&=&\sum_{w \in \mathcal{M}_0(m,n-m+1)}q^{{{\rm{maj}}}(w)}+\sum_{w \in \overline{\mathcal{M}}_1(m,n-m+1)}q^{{{\rm{maj}}}(w)}\\[5pt]
&=&q^m{n\brack m}+{n-1\brack m-1}\\[5pt]
&=&\frac{1-q^{n+m}}{1-q^n}{n\brack m}\\[5pt]
&=&{K}_q({n,m}).
\end{eqnarray*}
Furthermore, it is easy to check that \eqref{main-CombInt-e1} and \eqref{main-CombInt-e2} are valid when $m=0$. This completes the proof of Theorem \ref{main-CombInt}. \qed

\section{Proof of Theorem  \ref{main-q-log-cancave}}

Before we prove Theorem \ref{main-q-log-cancave}, it is useful to preset the following result.

\begin{lem}\label{Bulter-g}For $1\leq m\leq l< N$ and $M-m\geq N-l\geq 1,$
\[D_q(M,N,m,l)={M\brack m}{N\brack l}-
{M\brack m-1}{N\brack l+1}\]
has nonnegative coefficients as a polynomial in $q$.
\end{lem}

Lemma \ref{Bulter-g} reduces to the strong $q$-log-concavity of     Gaussian polynomials  when  $M=N$. We prove Lemma \ref{Bulter-g} by generalizing Butler's bijection \cite{Butler-1990}. To describe the proof, we need to recall   some  notation and terminology on partitions as  in \cite[Chapter 1]{Andrews-1976}.   {
A partition} $\lambda$ of a positive integer $n$ is a finite
nonincreasing sequence of positive integers
$(\lambda_1,\,\lambda_2,\ldots,\,\lambda_r)$ such that
$\sum_{i=1}^r\lambda_i=n.$  Then $\lambda_i$ are called the parts
of $\lambda$ and $\lambda_1$ is its largest part. The number of parts
of  $\lambda$ is called the length of $\lambda$, denoted by
$l(\lambda).$  The weight of $\lambda$ is the sum of  parts of $\lambda$, denoted
 $|\lambda|.$  The conjugate $\lambda'=(\lambda_1',\lambda_2',\ldots,\,\lambda'_t)$ of a partition $\lambda$ is defined by setting $\lambda'_i$ to be the number of parts of
$\lambda$ that are greater than or equal to $i$. Clearly,
$l(\lambda)=\lambda'_1$ and $\lambda_1=l(\lambda')$.

Let  $\mathcal{P}(m,n-m)$  denote  the set of partitions $\lambda$ such that $\ell(\lambda)\leq m$ and $\lambda_1\leq n-m$. It is well-known that the Gaussian polynomial  has the following partition interpretation \cite[Theorem 3.1]{Andrews-1976}:
\begin{equation}\label{int-GassCoef}
{n \brack m}=\sum_{\lambda \in \mathcal{P}(m,n-m)} q^{|\lambda|}.
\end{equation}

We are now prepared for the proof of Lemma \ref{Bulter-g} based on \eqref{int-GassCoef}.

\noindent{\it Proof of Lemma \ref{Bulter-g}:} For $1\leq m\leq l<N$ and $M-m\geq N-l\geq 1,$ by \eqref{int-GassCoef},  it suffices to construct an injection ${\Phi}$ from $\mathcal{P} ({m-1,M-m+1}) \times \mathcal{P} ({l+1,N-l-1}) $ to $\mathcal{P} ({m,M-m}) \times \mathcal{P} ({l,N-l}) $ such that if $\Phi(\lambda,\mu)=(\eta,\rho)$, then $|\lambda|+|\mu|=|\eta|+|\rho|$.

Let
\[\lambda=(\lambda_1,\lambda_2,\ldots,\lambda_{m-1})\in \mathcal{P} ({m-1,M-m+1}) \]
 and
 \[\mu=(\mu_1,\mu_2,\ldots,\mu_{l+1})\in \mathcal{P} ({l+1,N-l-1}),\]
 where $\lambda_1\leq M-m+1$ and $\mu_1\leq N-l-1$.

 We aim to construct a pair of partitions
  \[(\eta,\rho)\in \mathcal{P} ({m,M-m})\times \mathcal{P} ({l,N-l}).\]
   Let $I$ be the largest integer such that $\lambda_I\geq \mu_{I+1}+l-m+M-N+1$. If no such $I$ exists, then let $I=0$. In this case, we see that  $\lambda_1<M-m$ and set $\gamma=\lambda$ and $\tau=\mu$. Obviously,  $\gamma_1< M-m$ and $\tau_1<N-l$. We now assume that  $1\leq I\leq m-1$ and  define
\begin{equation}\label{defi-gamm}
\gamma=(\mu_1+(l-m+M-N+1), \ldots,\mu_I+(l-m+M-N+1),\lambda_{I+1},\ldots, \lambda_{m-1})
\end{equation}
and
\begin{equation}\label{defi-tau}
\tau=(\lambda_1-(l-m+M-N+1),\ldots,\lambda_I-(l-m+M-N+1),\mu_{I+1},\ldots, \mu_{l+1}).
\end{equation}
Since $I$ is the largest integer such that  $\lambda_I\geq \mu_{I+1}+(l-m+M-N+1)$, we get
\[\lambda_{I+1}< \mu_{I+2}+(l-m+M-N+1)\leq \mu_{I}+(l-m+M-N+1).\] It follows that $\gamma$ defined in \eqref{defi-gamm} and  $\tau$ defined in \eqref{defi-tau} are partitions. Furthermore,
\[\gamma_1=\mu_1+(l-m+M-N+1)\leq M-m \]
 {and}
\[    \tau_1=\lambda_1-(l-m+M-N+1)\leq N-l.\]

Let $\gamma'$ and $\tau'$ be the conjugates of $\gamma$ and $\tau$, respectively. We see that
\[\ell(\gamma')=\gamma_1\leq M-m \quad \text{and} \quad \ell(\tau')=\tau_1\leq N-l,\]
 so we can assume that
 \[\gamma'=(\gamma'_1,\gamma'_2,\ldots, \gamma'_{M-m})\]
  and
  \[\tau'=(\tau'_1,\tau'_2,\ldots, \tau'_{N-l}).\]
   Then
   \[\gamma'_1\leq m-1 \quad  \text{and} \quad \tau_1'\leq l+1.\]
     Let $J$ be the largest integer such that $\tau'_J\geq \gamma'_{J+1}+l-m+1$. If no such $J$ exists, let $J=0$, then $\tau'_1<l$ and set $\widetilde{\gamma}=\gamma',$ and $\widetilde{\tau}=\tau'$. Obviously,  $\widetilde{\gamma}_1< m$ and $\widetilde{\tau}_1<l$. We now assume that  $1\leq J\leq N-l$ and  define
\begin{equation}\label{conGam}
\widetilde{\gamma}=(\tau'_1-(l-m+1),\tau'_2-(l-m+1),\ldots,\tau'_J-(l-m+1),\gamma'_{J+1},\ldots, \gamma'_{M-m})
\end{equation}
and
\begin{equation}\label{contau}
\widetilde{\tau}=(\gamma'_1+(l-m+1),\gamma'_2+(l-m+1),\ldots,\gamma'_J+(l-m+1),\tau'_{J+1},\ldots, \tau'_{N-l}).
\end{equation}
Similarly, since $J$ is the largest integer such that  $\tau'_J\geq \gamma'_{J+1}+l-m+1$, we find that
\[\tau'_{J+1}<\gamma'_{J+2}+l-m+1\leq \gamma'_{J}+l-m+1,\]
 so $\widetilde{\gamma}$ defined in \eqref{conGam}  and $\widetilde{\tau}$ defined in \eqref{contau} are partitions. By the constructions of $\widetilde{\gamma}$ and $\widetilde{\tau}$, we see that \[\widetilde{\gamma}_1=\tau'_1-(l-m+1)\leq m\]
  and
  \[\widetilde{\tau}_1=\gamma'_1+(l-m+1)\leq l.\]
    Let $\eta$ and $\rho$ be the conjugates of $\widetilde{\gamma}$ and $\widetilde{\tau}$,  respectively. It is easy to check that $\eta \in \mathcal{P} ({m,M-m})$ and $\rho \in \mathcal{P} ({l,N-l})$. Furthermore,  this process is reversible. Thus, we complete the proof of Lemma \ref{Bulter-g}.    \qed

 Combining Lemma \ref{Bulter-g} and  the unimodality of  Gaussian polynomials, we  obtain the following   result.

\begin{lem}\label{Bulter-g2}For $1\leq m\leq l<N$, $M-m\geq N-l\geq 1$ and $0\leq r\leq M-N+2l-2m+2$,
\begin{equation}\label{diff-g}
D^r_q(M,N,m,l)={M\brack m}{N\brack l}-q^r
{M\brack m-1}{N\brack l+1}
\end{equation}
has nonnegative coefficients as a polynomial in $q$.
\end{lem}

\pf Let $A$ denote the degree of the polynomial ${M\brack m}{N\brack l}$ and let $B$ denote the degree of the polynomial ${M\brack m-1}{N\brack l+1}$. We have
\[A=m(M-m)+l(N-l), \]
\[ B=(m-1)(M-m+1)+(l+1)(N-l-1).\]
Furthermore, 
 \[A-B=M-N+2l-2m+2.\]
Let
\[{M\brack m}{N\brack l}=\sum_{i=0}^{A}a_i q^i, \quad {M\brack m-1}{N\brack l+1}=\sum_{i=0}^{B}b_i q^i\]
and let
\[D^r_q(M,N,m,l)={M\brack m}{N\brack l}-q^r{M\brack m-1}{N\brack l+1}=\sum_{i=0}^{A}c_iq^i,\]
where $c_i=a_i$ for $0\leq i< r$, $c_i=a_i-b_{i-r}$ for $r\leq i\leq B+r$ and $c_i=a_i$ for $B+r+1\leq i\leq A$.  It is easy to see that $c_i\geq 0$ for $0\leq i<r$ and  $B+r+1\leq i\leq A$. It remains to show that $c_i\geq 0$ for $r\leq i\leq B+r$.

It is known that the Gaussian polynomial ${M\brack m}$ is  symmetric and unimodal, see, for example, \cite[Theorem 3.10]{Andrews-1976} and \cite[Exercise 7.75]{Stanley-2012}, so
 \begin{equation}\label{lemsya}
 a_i=a_{A-i} \, \text{ for } \, 0\leq i\leq A, \quad \text{and} \quad  b_i=b_{B-i} \, \text{ for } \, 0\leq i\leq B,
 \end{equation}
\begin{equation} \label{lemunia}
a_0\leq a_1\leq \cdots \leq a_{\lfloor A/2\rfloor}=a_{\lceil A/2\rceil}\geq \cdots \geq a_{A-1}\geq a_{A},
\end{equation}

and
\begin{equation} \label{lemunib}
b_0\leq b_1\leq \cdots \leq b_{\lfloor A/2\rfloor}=b_{\lceil A/2\rceil} \geq \cdots \geq b_{B-1}\geq b_{B}.
\end{equation}

By Lemma \ref{Bulter-g}, we see that for $0\leq i\leq A$,
\begin{equation} \label{lemab}
a_i-b_i\geq 0.
\end{equation}
We consider the following  two cases:

Case 1. If $r\leq i\leq A/2$, then
\[c_i=a_i-b_{i-r}=a_i-a_{i-r}+a_{i-r}-b_{i-r},\]
which is nonnegative by \eqref{lemunia} and \eqref{lemab}.

Case 2. If $A/2\leq i\leq B+r$, then
\[c_i=a_i-b_{i-r}\overset{\eqref{lemsya}}{=}
a_{A-i}-b_{B-i+r}=a_{A-i}-a_{B-i+r}+a_{B-i+r}-b_{B-i+r},\]
which is nonnegative by \eqref{lemunia} and \eqref{lemab}.
Thus, we complete the proof of Lemma \ref{Bulter-g2}. \qed

We conclude this paper with a proof of Theorem \ref{main-q-log-cancave} by using Lemma \ref{Bulter-g2}.

\noindent{\it Proof of Theorem \ref{main-q-log-cancave}:} Recall that
\begin{align*}
K_q(n,m)&=\frac{1-q^{n+m}}{1-q^{n}}{n \brack m}={n \brack m}+q^n{n-1\brack m-1}.
\end{align*}
Hence
\begin{align*}
&K_q(n,m)K_q(n,l)-q^rK_q(n,m-1)K_q(n,l+1)\\[5pt]
&\quad =\left({n \brack m}+q^n{n-1\brack m-1}\right)\left({n \brack l}+q^n{n-1\brack l-1}\right)\\[5pt]
&\quad \quad -q^r\left({n \brack m-1}+q^n{n-1\brack m-2}\right)\left({n \brack l+1}+q^n{n-1\brack l}\right)\\[5pt]
&\quad = {n \brack m}{n \brack l}-q^r{n \brack m-1}{n \brack l+1} \\[6pt]
&\quad \quad +q^n\left( {n-1 \brack m-1}{n \brack l}-q^r{n-1 \brack m-2}{n \brack l+1} \right)\\[5pt]
&\quad \quad +q^n\left( {n \brack m}{n-1 \brack l-1}-q^r{n \brack m-1}{n-1 \brack l} \right)\\[5pt]
&\quad \quad +q^{2n}\left( {n-1 \brack m-1}{n-1 \brack l-1}-q^r{n-1 \brack m-2}{n-1 \brack l} \right).
\end{align*}
Using the notation in Lemma \ref{Bulter-g2}, we see that
\begin{align*}
&K_q(n,m)K_q(n,l)-q^rK_q(n,m-1)K_q(n,l+1)\\[5pt]
&=D^r_q(n,n,m,l)+q^nD^r_q(n-1,n,m-1,l)+q^nD^r_q(n,n-1,m,l-1)\\[5pt]
&\quad \quad +q^{2n}D^r_q(n-1,n-1,m-1,l-1).
\end{align*}
Applying Lemma 3.2, we find that  for $1\leq m\leq l<n$ and $0\leq r\leq 2l-2m+2$, 
\[D^r_q(n,n,m,l),\,D^r_q(n-1,n,m-1,l),\, \text{and}\  D^r_q(n-1,n-1,m-1,l-1) \]
have nonnegative coefficients as  polynomials in $q$, respectively, and  for $1\leq m\leq l<n$ and $0\leq r\leq 2l-2m+1$,
\[D^r_q(n,n-1,m,l-1)\]
has nonnegative coefficients as a polynomial in $q$. It follows that for $1\leq m\leq l<n$ and $0\leq r\leq 2l-2m+1$,
\begin{equation}\label{dif-K}
K_q(n,m)K_q(n,l)-q^rK_q(n,m-1)K_q(n,l+1)
\end{equation}
 has nonnegative coefficients as a polynomial in $q$. Hence it remains to show that the difference \eqref{dif-K}  has nonnegative coefficients as a polynomial in $q$ when $r=2l-2m+2$. It suffices to show that
\begin{equation}\label{uucc}
q^nD^{2l-2m+2}_q(n-1,n-1,m-1,l-1)+D^{2l-2m+2}_q(n,n-1,m,l-1)
\end{equation}
has nonnegative coefficients as a polynomial in $q$. First, it is easy to check that
\begin{eqnarray*}\label{middle}
&&q^nD^{2l-2m+2}_q(n-1,n-1,m-1,l-1)+D^{2l-2m+2}_q(n,n-1,m,l-1)\nonumber\\[8pt]
&&\quad=K_q(n,m){n-1\brack l-1}-q^{2l-2m+2}K_q(n,m-1){n-1\brack l}.
 \end{eqnarray*}
 Using the following relation:
 \begin{align*}
K_q(n,m)&=\frac{1-q^{n+m}}{1-q^{n}}{n \brack m}={n-1 \brack m-1}+q^m{n\brack m},
\end{align*}
we find that
\begin{eqnarray*}\label{middle}
&&q^nD^{2l-2m+2}_q(n-1,n-1,m-1,l-1)+D^{2l-2m+2}_q(n,n-1,m,l-1)\nonumber\\[8pt]
&&\quad=\left({n-1 \brack m-1}+q^m{n\brack m}\right){n-1\brack l-1}-q^{2l-2m+2}\left({n-1 \brack m-2}+q^{m-1}{n\brack m-1}\right){n-1\brack l}\nonumber\\[8pt]
&&\quad={n-1 \brack m-1}{n-1\brack l-1}-q^{2l-2m+2}{n-1 \brack m-2}{n-1\brack l}\nonumber\\[8pt]
&&\quad \quad \quad +q^{m}\left({n\brack m}{n-1\brack l-1}-q^{2l-2m+1}{n\brack m-1}{n-1\brack l}\right)\nonumber\\[8pt]
&&\quad =D^{2l-2m+2}_q(n-1,n-1,m-1,l-1)+q^mD^{2l-2m+1}_q(n,n-1,m,l-1).
 \end{eqnarray*}
 From Lemma \ref{Bulter-g2}, we see that
 \[D^{2l-2m+2}_q(n-1,n-1,m-1,l-1), \quad \text{and} \quad D^{2l-2m+1}_q(n,n-1,m,l-1)\]
  have nonnegative coefficients as a polynomial in $q$, respectively, and so
  \eqref{uucc} has nonnegative coefficients as a polynomial in $q$.  Thus, we complete the proof of Theorem \ref{main-q-log-cancave}. \qed

 \vskip 0.5cm

 \noindent{\bf Acknowledgment.}  We wish to thank the referees for   invaluable comments and suggestions.  This work was supported by the National Science Foundation of
China.

\end{document}